\documentclass[reqno,12pt]{amsart}
\usepackage{amscd,amssymb}

\theoremstyle{plain}
\newtheorem{Thm}[subsection]{Theorem}
\newtheorem{Cor}[subsection]{Corollary}
\newtheorem{Lem}[subsection]{Lemma}
\newtheorem{Prop}[subsection]{Proposition}
\newtheorem{Conj}[subsection]{Conjecture}

\theoremstyle{definition}
\newtheorem{Def}[subsection]{Definition}

\theoremstyle{remark}

\newtheorem{Rem}[subsection]{Remark}

\newtheorem{conj}{Conjecture}

\numberwithin{equation}{section}
\renewcommand{\rm}{\normalshape}

\newif\ifShowLabels
\ShowLabelstrue
\newdimen\theight
\def\TeXref#1{%
        \leavevmode\vadjust{\setbox0=\hbox{{\tt
                \quad\quad  {\small \rm #1}}}%
        \theight=\ht0
        \advance\theight by \lineskip
        \kern -\theight \vbox to
        \theight{\rightline{\rlap{\box0}}%
        \vss}%
        }}%

\ShowLabelsfalse

\renewcommand{\sec}[2]{\section{#2}\label{S:#1}%
        \ifShowLabels \TeXref{{S:#1}} \fi}
\newcommand{\ssec}[2]{\subsection{#2}\label{SS:#1}%
        \ifShowLabels \TeXref{{SS:#1}} \fi}

\newcommand{\refss}[1]{Section ~\ref{SS:#1}}

\newcommand{\reft}[1]{Theorem ~\ref{T:#1}}
\newcommand{\refl}[1]{Lemma ~\ref{L:#1}}
\newcommand{\refp}[1]{Proposition ~\ref{P:#1}}

\newcommand{\refe}[1]{\eqref{E:#1}}
\newcommand{\refco}[1]{Conjecture ~\ref{Co:#1}}

\newenvironment{thm}[1]%
        { \begin{Thm} \label{T:#1}  \ifShowLabels \TeXref{T:#1} \fi }%
        { \end{Thm} }

\renewcommand{\th}[1]{\begin{thm}{#1} \sl }
\renewcommand{\eth}{\end{thm} }

\newenvironment{lemma}[1]%
        { \begin{Lem} \label{L:#1}  \ifShowLabels \TeXref{L:#1} \fi }%
        { \end{Lem} }
\newcommand{\lem}[1]{\begin{lemma}{#1} \sl}
\newcommand{\elem}{\end{lemma}}

\newenvironment{propos}[1]%
        { \begin{Prop} \label{P:#1}  \ifShowLabels \TeXref{P:#1} \fi }%
        { \end{Prop} }
\newcommand{\prop}[1]{\begin{propos}{#1}\sl }
\newcommand{\eprop}{\end{propos}}

\newenvironment{corol}[1]%
        { \begin{Cor} \label{C:#1}  \ifShowLabels \TeXref{C:#1} \fi }%
        { \end{Cor} }
\newcommand{\cor}[1]{\begin{corol}{#1} \sl }
\newcommand{\ecor}{\end{corol}}

\newenvironment{defeni}[1]%
        { \begin{Def} \label{D:#1}  \ifShowLabels \TeXref{D:#1} \fi }%
        { \end{Def} }
\newcommand{\defe}[1]{\begin{defeni}{#1} \sl }
\newcommand{\edefe}{\end{defeni}}

\newenvironment{remark}[1]%
        { \begin{Rem} \label{R:#1}  \ifShowLabels \TeXref{R:#1} \fi }%
        { \end{Rem} }
\newcommand{\rem}[1]{\begin{remark}{#1}}
\newcommand{\erem}{\end{remark}}

\newenvironment{conjec}[1]%
        { \begin{Conj} \label{Co:#1}  \ifShowLabels \TeXref{Co:#1} \fi }%
        { \end{Conj} }
\renewcommand{\conj}[1]{\begin{conjec}{#1} \sl }
\newcommand{\econj}{\end{conjec}}

\newcommand{\eq}[1]%
        { \ifShowLabels \TeXref{E:#1} \fi
           \begin{equation} \label{E:#1} }
\newcommand{\eeq}{ \end{equation} }

\newcommand{\prf}{ \begin{proof} }
\newcommand{\epr}{ \end{proof} }

\newcommand\alp{\alpha}         

\newcommand\gam{\gamma}         
\newcommand\del{\delta}

\newcommand\calB{{\mathcal{B}}}

\newcommand\calD{{\mathcal{D}}}
\newcommand\calE{{\mathcal{E}}}
\newcommand\calF{{\mathcal{F}}}
\newcommand\calG{{\mathcal{G}}}
\newcommand\calH{{\mathcal{H}}}

\newcommand\calK{{\mathcal{K}}}

\newcommand\calM{{\mathcal{M}}}
\newcommand\calN{{\mathcal{N}}}

\newcommand\QQ{\mathbb{Q}}

\renewcommand\AA{\mathbb{A}}

\newcommand\sdp{\times \hskip -0.3em {\raise 0.3ex
\hbox{$\scriptscriptstyle |$}}} 

\newcommand\Perv{\operatorname{Perv}}

\newcommand\tilG{{\widetilde{G}}}

\newcommand\x{\times}
\newcommand\ten{\otimes}

\newcommand\qlb{{\overline \QQ}_l}

\begin{document}
\title{On Euler characteristic of equivariant sheaves}
\author{Alexander Braverman}
\address{Department of Mathematics, Harvard University, 1 Oxford st.
Cambridge MA, 02138}
\email{braval\@math.harvard.edu}
\begin{abstract}
Let $k$ be an algebraically closed field of characteristic $p>0$ and let
$\ell$ be another prime number. O.~Gabber and F.~Loeser proved that for any
algebraic torus $T$ over $k$ and any perverse $\ell$-adic sheaf $\calF$
on $T$ the Euler characteristic $\chi(\calF)$ is non-negative.

We conjecture that the same result holds for any perverse sheaf $\calF$ on a
reductive group $G$ over $k$ which is equivariant with respect to the
adjoint action. We prove the conjecture when $\calF$ is obtained by Goresky-MacPherson extension
from the set of regular semi-simple elements in $G$. From this we deduce that the
conjecture holds  for $G$ of semi-simple rank 1.
\end{abstract}
\maketitle

\sec{}{The main conjecture}

\ssec{}{Notations}
In what follows $k$ denotes an algebraically closed field of characteristic
$p>0$. Let $\ell$ be a prime number different from $p$. For an algebraic
variety $X$ over $k$ we denote by $\calD(X)$ the derived category of
$\ell$-adic sheaves on $X$. Also we denote by $\Perv(X)\subset \calD(X)$ the
subcategory of perverse sheaves. For any $\calF\in\calD(X)$ we denote by
$\chi(\calF)$ the  Euler characteristic of $\calF$, i.e.
\eq{}
\chi(\calF)=\sum (-1)^i\dim H^i_c(X,\calF)=\sum (-1)^i\dim H^i(X,\calF)
\end{equation}
(cf. \cite{La} for the proof of the equality).

Let $T$ be an algebraic torus over $k$. The following theorem is proved in
\cite{GL} (cf. also \cite{FK} for a partial analogue in characteristic zero):
\th{torus}
Let $\calF\in\Perv(T)$. Then $\chi(\calF)\geq 0$.
\eth
Let $G$ be a connected reductive algebraic group over $k$. We shall denote by
$G_{rs}\subset G$ the open subspace of regular semi-simple elements.
We denote by $\Perv_G(G)$ the subcategory of $\Perv(G)$ consisting of
perverse sheaves which are equivariant with respect to the adjoint action.
We propose the following generalization of \reft{torus}.
\conj{main}
Let $\calF\in\Perv_G(G)$. Then $\chi(\calF)\geq 0$.
\econj
\th{main}
Assume that $\calF\in\Perv_G(G)$ is equal to the Goresky-MacPherson
extension of its restriction to $G_{rs}$. Then $\chi(\calF)\geq 0$.
\eth
\cor{}
\refco{main} holds for $G$ of semi-simple rank 1.
\ecor
\prf
Let $\calN$ be the cone of unipotent elements in $G$ and let $Z$ be the
center of $G$. Clearly it is enough to prove \refco{main} only for irreducible sheaves.
However if $\calF\in\Perv_G(G)$ is irreducible then either $\calF$ is equal to the
Goresky-MacPherson extension of its restriction to $G_{rs}$ or it is supported on
$Z\calN$. The former case is covered by \reft{main}; hence it is enough to deal with the
latter.

First of all there exists a connected component $Z'$ of $Z$ such that
$\calF$ is supported on $Z'\calN$. Thus  one of the following
is true:

1) $\calF$ is supported on $Z'$.

2) $\calF$ is equal to $\calF'\boxtimes (\qlb)_\calN[2]$ where $\calF$ is
an irreducible perverse sheaf on $Z'$ and $(\qlb)_{\calN}$ is the constant
sheaf on $\calN$.

3) $\calF=\calF'\boxtimes \calE$ where $\calF'$ is an irreducible perverse
sheaf on $Z'$ and $\calE$ is an irreducible perverse sheaf on $\calN$ whose
restriction to the open orbit is isomorphic to (the only) non-trivial
equivariant irreducible local system on this orbit (this local system is of rank 1 and
its square is isomorphic to the constant sheaf).

In cases 1 and 2 our result follows immediately from \reft{torus}
(note that $Z'$ is a torus and that in case 2 we have
$\chi(\calF)=\chi(\calF')$).
In case 3 it is known that $H^*(\calN,\calE)=0$, hence
$H^*(G,\calF)=0$, hence $\chi(\calF)=0$.
\epr

\noindent
{\it Remark.} \reft{torus} has an "ideological"
explanation: namely in \cite{GL} O.~Gabber and F.~Loeser construct a
"Mellin transform" functor $\calM:\calD^b(T)\to \calD^b_{coh}(\text{Loc}_T)$ where
$\text{Loc}_T$ denotes the moduli space of tame rank one local systems on $T$ and
$\calD^b_{coh}(\text{Loc}_T)$ is the derived category of quasi-coherent sheaves
on $\text{Loc}_T$. Moreover, for every $\calF\in\Perv(T)$ the complex
$\calM(\calF)$ is actually a sheaf and $\chi(\calF)$ is equal to the generic
rank of $\calM(\calF)$ (hence $\chi(\calF)\geq 0$). We don't know whether a
similar explanation of \refco{main} is possible.

We expect that \refco{main} holds also when $k$ has characteristic 0 and perverse
sheaves are replaced by  holonomic D-modules. When $G$ is a torus and the
D-modules in question have regular singularities a beautiful geometric proof
of \reft{torus} was given in \cite{FK} where the authors explain how to
compute the corresponding Euler characteristics using certain non-compact
generalization of Kashiwara's index theorem. It would be interesting to
generalize the proof of \cite{FK} to the case of arbitrary reductive group. Some results in this
direction have recently been obtained by V.Kiritchenko.
\sec{}{Proof of \reft{main}}

\ssec{}{The space $\tilG$}Let $\calB$ denote the flag variety of $G$, i.e. the variety of all Borel
subgroups of $G$. Let $\tilG$ denote the variety of pairs $(B\in\calB,g\in
B)$. We have the natural maps $\pi:\tilG\to G$, $\alp:\tilG\to T$ and
$\beta:\tilG\to\calB$. It is easy to see that $\alp$ is smooth and $\pi$ is proper.
Let $d=\dim G-\dim T=2\dim\calB$. Set
$\tilG_{rs}=\pi^{-1}(G_{rs})=\alp^{-1}(T_{rs})$.

\lem{induction}
\begin{enumerate}
\item
Let $\calG\in\Perv(T)$. Assume that $\calG$ is equal to the Goresky-MacPherson extension
of its restriction
on $T_{rs}$.
$\pi_!\alp^*\calG[d](\frac{d}{2})$ is a perverse sheaf which is equal to the Goresky-MacPherson extension
of its restriction to $G_{rs}$.
\item
Every $W$-equivariant structure on a sheaf $\calG\in\Perv(T)$ as above gives rise to
a $W$-action on $\pi_!\alp^*\calG[d](\frac{d}{2})$.
\item
Let $\calF\in\Perv_G(G)$. Assume that $\calF$ is equal to the
Goresky-MacPherson extension of its restriction to $G_{rs}$.
Then there exists a $W$-equivariant sheaf $\calG\in\Perv(T)$ which is equal to the
Goresky-MacPherson extension of its restriction to $T_{rs}$ such that
\eq{fginvariant}
\calF=(\pi_!\alp^*\calG[d](\frac{d}{2}))^W.
\end{equation}
\end{enumerate}
\elem
\prf The first statement of \refl{induction} is well-known (it follows from
the smallness property of $\alp$ -- cf. \cite{CG} and references therein)

It follows from 1 that in order to prove 2 it is enough to construct an action of $W$ on
the restriction of $\pi_!\alp^*\calG[d](\frac{d}{2})$ to $G_{rs}$. Note that $W$ acts
freely on $\tilG_{rs}$ in the natural way and
\eq{quotient}
G_{rs}=\tilG_{rs}/W.
\end{equation}
Moreover, the restriction of $\alp$ to $\tilG_{rs}$ is $W$-equivariant.
Hence $\alp^*\calG[d](\frac{d}{2})|_{\tilG_{rs}}$ is $W$-equivariant and thus
\refe{quotient} implies that $\pi_!\alp^*\calG[d](\frac{d}{2})|_{G_{rs}}$ has
a natural action of $W$.

Let us prove 3. Choose an embedding $i:T\to G$. Let $i_{rs}:T_{rs}\to G$ be
its restriction to $T_{rs}$. Let $\calG_{rs}=i_{rs}^*\calF[-d](-\frac{d}{2})$. It follows from
the $G$-equivariance of $\calF$ that $\calG_{rs}$ is perverse. We let
$\calG$ be its Goresky-MacPherson extension to $T$. We leave the
verification of \refe{fginvariant} to the reader.
\epr
\prop{torus-g}
Let $\calF$ be as in
\reft{main} and let $\calG$ be as in \refl{induction}(3). Then
$\chi(\calF)=\chi(\calG)$.
\eprop
Clearly \refp{torus-g} together with \reft{torus} imply \reft{main}.
\ssec{proof}{Proof of Proposition \refp{torus-g}}
Set $\calH=\alp^*\calG[d](\frac{d}{2})\in\Perv(\tilG)$. Clearly,
$H^*_c(G,\calF)=H^*_c(\tilG,\calH)^W$.

Consider $\alp_!(\calH)$. By the projection formula it is isomorphic
to $\calG\ten \alp_!(\qlb)_{\tilG}[d](\frac{d}{2})$. It is easy to see that the
complex $\alp_!(\qlb)_{\tilG}[d](\frac{d}{2})$ is
constant with fibers isomorphic to $H^*(\calB,\qlb)$. Indeed, let $\del:\calB\x T\to T$
denote the projection to the first multiple.
Then $\alp_!(\qlb)_{\tilG}=\del_!(\beta\x\alp)_!(\qlb)_{\tilG}$. However,
the map $(\beta\x \alp):\tilG\to\calB\x T$ is a locally trivial fibration
(in Zariski topology) with fiber isomorphic to $\AA^d$. Thus
$(\beta\x \alp)_!(\qlb)_{\tilG}\simeq (\qlb)_{\calB\x T}[-d](-\frac{d}{2})$ and
hence
\eq{projection}
\alp_!(\qlb)_{\tilG}\simeq (\qlb)_T\ten H^*(\calB,\qlb)[-d](-\frac{d}{2}).
\end{equation}
This clearly gives rise to the isomorphism
\eq{twosides}
H^*_c(\tilG,\calH)\simeq H^*_c(T,\calG)\ten H^*(\calB,\qlb).
\end{equation}
We can also characterize the isomorphism \refe{projection}
in the following way. To determine it uniquely it is enough to construct it on $T_{rs}$.
Choose as before an embedding
$i:T\to G$. There is a canonical Borel subgroup $B$ containing $i(T)$ (recall that the abstract
Cartan group $T$ comes with a root system with a preferred set of positive roots). Then for every
$t\in T_{rs}$ we have canonical isomorphism
\eq{}
\alp^{-1}(t)\simeq G/T
\end{equation}
(which depends however on the above choice). Clearly
$H^*_c(G/T,\qlb)[d](\frac{d}{2})=H^*(\calB,\qlb)$. Hence we have constructed an isomorphism
between $\alp_!(\qlb)_{\tilG}[d](\frac{d}{2})|_{T_{rs}}$ and the constant complex with
fiber $H^*(\calB,\qlb)$. It is easy to see that this isomorphism does not
depend on the choice of the embedding $T\to G$ and coincides with
\refe{projection}.

Let $W$ act on the left hand side by means
of the identification $H^*_c(\tilG,\calH)=H^*_c(G,\pi_!\calH)$ (note that
by \refl{induction} the group $W$ acts on $\pi_!\calH$).
Let $W$ also act on the right hand side by means of the tensor product of
the $W$-action on $H^*_c(T,\calG)$ coming from the $W$-equivariant structure
on $\calG$ and the natural $W$-action on $H^*(\calB,\qlb)$.
\lem{w}The isomorphism \refe{twosides}
is also an isomorphism of $W$-modules with respect to the above actions.
\elem
\prf
Let $Z=G\underset{T/W}\x T$ be the image of $\tilG$ in $G\x T$ under the
map $\gam=\pi\x \alp$. This is a closed subvariety of $G\x T$.
It is invariant with respect to the $W$ action on the second multiple.
Let $Z_{rs}$ denote the set of "regular semi-simple" elements of $Z$
(i.e. the set of all elements of $Z$ whose projection to $Z$ is regular
semi-simple).

Let $\calK=\gam_!(\calH)$. Then $\calK$ is equal to the Goresky-MacPherson extension
of its restriction to $Z_{rs}$. Since $\gam$ induces an isomorphism
$\tilG_{rs}\simeq Z_{rs}$ and since the restriction of $\calH$ to $\tilG_{rs}$ is
$W$-equivariant it follows that $\calK$ also has a natural $W$-equivariant structure
as a perverse sheaf on $G\x T$ (where the $W$-action is as before on the
second multiple). Thus $W$ acts naturally on $H^*_c(G\x T,\calK)$.
We have the natural isomorphism $H^*_c(\tilG,\calH)\simeq
H^*_c(G\x T,\calK)$ and by the definition the action of $W$ on
$H^*_c(\tilG,\calH)$ introduced before \refl{w} corresponds to the
action of $W$ on $H^*_c(G\x T,\calK)$ introduced above.

Consider, on the other hand, the complex $p_!(\calK)$ where $p:G\x T\to T$
denotes the natural projection. Clearly, we have $p_!(\calK)=\alp_!(\calH)\simeq
\calG\ten H^*(\calB,\qlb)$. On the other hand, since $p$
commutes with $W$ it follows that $p_!(\calK)$ admits a natural
$W$-equivariant structure. To prove \refl{w} it is enough to show that
under the identification $p_!(\calK)\simeq \calG\ten
H^*(\calB,\qlb)[d](\frac{d}{2})$ this structure is equal to the tensor
product of the original $W$-equivariant structure on $\calG$ with the
natural $W$-action on $H^*(\calB,\qlb)$. Since every perverse cohomology
of $p_!(\calK)$ is equal to the Goresky-MacPherson extension of its
restriction to $T_{rs}$ it follows that it is enough to check this equality
only on $T_{rs}$. This, however, follows immediately from the description of
the $W$-equivariant structure on
$\alp_!(\qlb)_{\tilG}[d](\frac{d}{2})|_{\tilG_{rs}}$ given in \refss{proof}
and the following observation:

Choose as before an embedding $T\to G$.
Let $N(T)$ denote the normalizer of $T$ in $G$. Then $W=N(T)/T$ acts
on $G/T$ and hence on $H^*_c(G/T,\qlb)$. By identifying
$H^*_c(G/T,\qlb)[d](\frac{d}{2})$ with $H^*(\calB,\qlb)$ we get an
action of $W$ on the latter. This is the standard $W$-action on
$H^*(\calB,\qlb)$.
\epr

It follows from \refl{w} that $\chi(\calF)$ is equal to the Euler characteristic of
$(H^*_c(\calG)\ten H^*(\calB,\qlb))^W$.
However, it is known that $H^*(\calB,\qlb)$ lives only in even degrees and
when we forget the grading it is isomorphic to the regular representation $\qlb[W]$ of
$W$. Hence
\begin{align*}
\chi(\calF)=\chi((H^*_c(T,\calG)\ten& H^*(\calB,\qlb)[d])^W)=\\
&\chi((H^*_c(T,\calG)\ten \qlb[W])^W)=\chi(H^*_c(T,\calG))=\chi(\calG)
\end{align*}
which finishes the proof.
\ssec{}{Acknowledgements}
The author is grateful to F.~Loeser for an illuminating discussion on the subject
held in the summer of 2001, to V.~Kiritichenko for making her results available to the
author before publication and to M.~Kapranov for explaining the contents of
\cite{FK}.

\end{document}